\documentclass[12pt]{article} 
\usepackage[final]{graphicx}
\usepackage[dvips]{color}
\usepackage{graphicx}
\usepackage{amsmath}
\usepackage{amssymb}
\usepackage{ascmac}
\usepackage{amsthm}
\usepackage{bm}
\def\Bbb{\bf} 

\newcommand\C{{ \Bbb C}}

\newcommand\Z{{\Bbb Z}}
\newcommand\Q{{\Bbb Q}}

\newcommand\la{{\lambda}}

\newtheorem{lmm}{Lemma}[section]
\newtheorem{thm}[lmm]{Theorem}

\newtheorem{rmk}[lmm]{Remark}

\def\comment#1{ }

\makeatletter
    
    \@addtoreset{equation}{section}
  \makeatother
\allowdisplaybreaks
\begin{document}
\title{On a strange evaluation \\of the hypergeometric series by Gosper}
\author{Akihito Ebisu
}

\maketitle
\begin{abstract}
In a series of letters to D.Stanton, R.W.Gosper presented many strange evaluations of hypergeometric series.
Recently, we rediscovered one of the strange hypergeometric identities appearing in [Go].
In this paper, we prove this identity 
and derive its generalization using contiguity operators.


\end{abstract}
\tableofcontents
\section{Introduction and main theorem}
In a series of letters to D.Stanton, R.W.Gosper carried out many strange evaluations of hypergeometric series.
Some are presented in [GS].
However, the following identity appears in [Go], but not in [GS]:
\begin{gather}
F\left(1-a, b, b+2; \frac{b}{a+b}\right)=(b+1)\left(\frac{a}{a+b}\right)^a,
\end{gather} 
where 
$$
F(a,b,c;x)
:=\sum^{\infty}_{n=0}\frac{(a,n)(b,n)}{(c,n)(1,n)}\,x^n,
$$
and $(a,n):=\Gamma (a+n)/\Gamma (a)$.
Recently, we discovered (1.1) independently.
Here, we present our proof of this identity 
and derive a generalization using contiguity operators.

The main theorem in this article is the following:
\begin{thm}
We assume that $\ell \in \Z _{>0}, a \in \C$ and $c\in \C \setminus\Z$.
For any root $\la$ of $F(1-a,-\ell ,2-c;x)$,
which is a polynomial in the variable $x$ of degree at most $\ell$,
we have
\begin{align}
F(a, 1+\ell ,c;\la)&=-\frac{(1-c)q_0(\la)}{(1,\ell )(1-\la)^\ell },\\
F(c-a, c-1-\ell ,c; \la)&=-\frac{(1-c)}{(1,\ell )}(1-\la)^{a+1-c}q_0(\la),
\end{align}
where $q_0(x)$ is the polynomial in $x$ of degree at most $\ell -1$ given by
\begin{multline}
q_0(x)=-\frac{(1,\ell )}{1-c} (1-x)^{c-a-1}F(c-a, c-1-\ell ,c;x)\\
+\frac{(2-c,\ell )}{1-c} F(a, 1, c;x)F(1-a, -\ell ,2-c;x).
\end{multline}
For example, when $\ell =1$, we have
\begin{align}
F\left(a, 2, c;\frac{c-2}{a-1}\right)&=\frac{(a-1)(c-1)}{a+1-c},\\
F\left(c-a, c-2, c;\frac{c-2}{a-1}\right)&=(c-1)\left(\frac{a+1-c}{a-1}\right)^{a+1-c}.
\end{align}
\end{thm}
\begin{rmk}
If we assume further the condition $a \notin \Z$,
then the degree of $F(1-a,-\ell ,2-c;x)$ is exactly $\ell$. 
In this case, $q_0(x)$ can also be expressed as
\begin{multline}
q_0(x)=(2-a,\ell -1)(-x)^{\ell -1}\times \left(
{\text {the $\ell$-th partial sum of the power series in $1/x$}}\right. \\
\left. F(2-c, 1, 2-a; 1/x)F(c-1-\ell , -\ell , a-\ell ;1/x)\right),
\end{multline}
and thus, the degree of $q_0(x)$ is exactly $\ell-1$.
\end{rmk}
Note that (1.6) is equivalent to (1.1).

Many methods for discovering and proving hypergeometric identities are known.
In the 19th century, such identities were discovered using algebraic transformations of hypergeometric series
(cf. [V] for a treatment of algebraic transformations of hypergeometric series).
Moreover, during the last several decades, 
many new methods that exploit progress in computer technology have been constructed:
Gosper's algorithm, the W-Z method, Zeilberger's algorithm, etc. (cf. [Ko] and [PWZ]).
These well-known algorithms have been used for discovering and proving hypergeometric identities
expressed in closed forms (cf. [WoGo], [WoHy], [WoWZ] and [WoZe]).
Indeed, in [Ek] and [AZ], M.Apagodu, S.B.Ekhad and D.Zeilberger 
report the discovery of such identities using these algorithms.
In addition, such identities are proved with the aid of these algorithms in [Ko] and [PWZ].
However, note that although $q_0(x)$ appearing in (1.2) and (1.3) is expressed explicitly,
it is not in a closed form (see. (1.4)).
Therefore, 
it seems that identities (1.2) and (1.3) could not have been discovered and can not be proved
through a direct application of the methods mentioned above.
\begin{rmk}
If we input 
$$
{\tt >simplify(hypergeom([c-a, c-2], [c],(c-2)/(a-1)))}
$$
in Maple 16, we obtain (1.6) as an output.
From this, we can easily obtain (1.5) using Euler transformation
\begin{gather}
F(a,b,c;x)=(1-x)^{c-a-b}F(c-a,c-b,c;x)
\end{gather}
(cf. (2.2.7) in [AAR]).
However, when we input
$$
{\tt >simplify(hypergeom([a, 2], [c],(c-2)/(a-1)))},
$$
Maple 16 does not return (1.5).
This is mysterious.
\end{rmk}

\section{Contiguity operators}
In this section, we introduce contiguity operators 
and express a composition of these operators in terms of the hypergeometric operator. 
Modifying the results given in [Eb], which treats this problem, 
we obtain Lemmas 2.2 and 2.3.
These are used for proving the main theorem.

Let $\partial:=d/dx$.
Further, let $L(a, b, c)$ denote the hypergeometric (differential) operator in $x$
\begin{gather}
\partial ^2+\frac{c-(a+b+1)x}{x(1-x)}\partial -\frac{ab}{x(1-x)},
\end{gather}
where $a, b$ and $c$ are complex variables. 
We call the differential operator $(x\partial +b)$ a {\bf contiguity operator}
because the relation
\begin{gather}
(x\partial +b)F(a, b, c; x)=bF(a, b+1, c; x)
\end{gather}
holds for $c\notin \Z_{\leq 0}$ (cf. Proposition 2.1.2 in [IKSY]).

Now, we consider the composition of contiguity operators 
\begin{gather}
(x\partial +b+\ell -1)\cdots(x\partial +b+1)(x\partial +b)=:H(\ell ).
\end{gather}
Then, $H(\ell )$ can be expressed as
\begin{gather}
H(\ell)=p(\partial)L(a,b,c)+q(x)\partial +r(x),
\end{gather}
where $p(\partial)\in \Q (a, b, c, x)[\partial]$ and $q(x), r(x)\in \Q (a, b, c, x)$. 
Here, $\Q (a, b, c, x)$ is the rational function field in $a, b, c$ and $x$ over $\Q$,
and $\Q (a, b, c, x)[\partial]$ is the ring of the differential operators in $x$ over $\Q (a, b, c)$.

Further, we assume the following:
\begin{align*}
&{\text {A1:}}\ a, b, c-a, c-b\notin \Z  &{\text {(cf. Section 1 of [Eb])}},\\
&{\text {A2:}}\ c, c-a-b, a-b \notin \Z  &{\text {(cf. Section 1 of [Eb])}},\\
&{\text {E1:}}\ (b,\ell)-(b+1-c,\ell)\neq 0  &{\text {(cf. Lemma 3.1 of [Eb])}},\\
&{\text {E2':}}\ \ell \neq 1 \ {\rm or}\  \frac{(b+1,\ell -1)}{a}-\frac{(c-a,\ell -1)}{c-b-1}\neq 0  
&{\text {(cf. Subsection 3.2 of [Eb]).}}
\end{align*}
Assuming the above, we can directly apply propositions and theorems given in [Eb].
Specifically,  
we need only substitute $k=0$, $l=\ell$ and $m=0$ into the formulas
appearing in the propositions and theorems of [Eb].
Doing so, $q(x)$ and $r(x)$ can be expressed as
\begin{align}
\begin{split}
&q(x)=x^{v_0}(1-x)^{v_1}q_0(x),\\
&q_0(x):{\text{a polynomial in $x$ of degree $g$ and }} q_0(0)q_0(1)\neq 0,
\end{split}\\
\begin{split}
&r(x)=x^{w_0}(1-x)^{w_1}r_0(x),\\
&r_0(x):{\text{a polynomial in $x$ of degree $h$ and }} r_0(0)r_0(1)\neq 0
\end{split}
\end{align}
(cf. Section 3 of [Eb]). Moreover, by Propositions 3.4 and 3.9 in [Eb], we have
\begin{gather}
(v_0, v_1, g)=(1, 1-\ell , \ell -1),\ (w_0, w_1, h)=(0, 1-\ell , \ell -1),
\end{gather}
and by Theorems 3.7 and 3.10 in [Eb], we have
\begin{align}
q_0(x)&=-\frac{(b,\ell)}{1-c}F(c-a, c-b-\ell , c; x) F(a+1-c, b+1-c, 2-c; x)\notag \\
&\quad +\frac{(b+1-c,\ell )}{1-c}F(a, b, c; x)F(1-a, 1-b-\ell , 2-c; x) \\
r_0(x)&=(b,\ell)F(c-a, c-b-\ell , c;x)F(a+1-c, b+1-c, 1-c;x) \notag\\
&\quad -\frac{ab(b+1-c,\ell)}{c(1-c)}xF(a+1, b+1, c+1;x)F(1-a, 1-b-\ell , 2-c; x).
\end{align}
\begin{rmk}
Let us consider the case in which conditions A1, A2, E1 and E2' are not satisfied.
Even in this case, 
because $q(x)$ and $r(x)$ are rational functions in $a, b, c$ and $x$,
we have
\begin{gather}
q(x)=x(1-x)^{1-\ell}q_0(x),\quad
r(x)=(1-x)^{1-\ell}r_0(x),
\end{gather}
where $q_0(x)$ and $r_0(x)$ are polynomials in $x$ of degree at most $\ell -1$.
Moreover, if the right-hand sides of (2.8) and (2.9) are meaningful (that is, $c\notin \Z$),
then (2.8) and (2.9) still hold.
Note that $q_0(x)$ and $r_0(x)$ may vanish at $x=0$ or $1$ or may be of smaller degree than that stated in (2.7).
\end{rmk}
Now, we consider the case in which $a\in \C$, $b=1$ and $c\in \C \setminus\Z$.
Then, we obtain the following lemma from Remark 2.1:
\begin{lmm} 
Let $\ell$ be a positive integer.
We represent $(x\partial +\ell)\cdots(x\partial +2)(x\partial +1)$ by $H_1(\ell )$.
Then, expressing $H_1(\ell )$ in terms of the hypergeometric operator $L(a, 1, c)$,
we have 
\begin{gather}
H_1(\ell)=p(\partial)L(a, 1, c)+q(x)\partial +r(x).
\end{gather}
Here, $q(x)$ and $r(x)$ are given by
\begin{gather}
q(x)=x(1-x)^{1-\ell}q_0(x),\quad
r(x)=(1-x)^{1-\ell}r_0(x),
\end{gather}
where, $q_0(x)$ and $r_0(x)$ are polynomials of degree {\bf {at most}} $\ell -1$ given by 
\begin{align}
q_0(x)&=-\frac{(1,\ell)}{1-c} (1-x)^{c-a-1}F(c-a, c-1-\ell , c;x) \notag\\
&\quad +\frac{(2-c,\ell)}{1-c} F(a,1,c;x)F(1-a,-\ell ,2-c;x), \\
r_0(x)&=(1,\ell) F(c-a,c-1-\ell ,c;x)F(a+1-c, 2-c,1-c;x)\notag \\
&\quad -\frac{a(2-c,\ell )}{c(1-c)} xF(a+1, 2, c+1;x)F(1-a, -\ell ,2-c;x).
\end{align}
\end{lmm}
We close this section by presenting the following lemma:
\begin{lmm}
Writing
\begin{gather}
y_1(x):=F(a, 1, c; x),\ y_2(x):=x^{1-c}(1-x)^{c-a-1},
\end{gather}
we have
\begin{flalign}
&y_1(x)=-(1-c)y_2(x)\int_0^{x}\frac{1}{t(1-t)y_2(t)} dt \ \left(\Re c>1\right), \\
&H_1(\ell)y_1(x)=(1,\ell)F(a,1+\ell,c;x),\\
&H_1(\ell)y_2(x)=(2-c,\ell )y_2(x)(1-x)^{-\ell}F(1-a, -\ell ,2-c;x).
\end{flalign}
\end{lmm}
First, regarding (2.16), see (2.1) in [Du].   
Note that (2.16) implies that $F(a, 1, c; x)$ can be expressed in terms of the incomplete beta function.
Next, it is easily shown that (2.17) can be obtained from (2.2).
Finally, we now demonstrate (2.18). It is known that
\begin{multline}
(x\partial +b)x^{1-c}F(a+1-c, b+1-c, 2-c;x)\\
=(b+1-c)x^{1-c}F(a+1-c, b+2-c, 2-c;x)
\end{multline} 
(cf. Lemma 2.2 in [Eb]).
Then, applying the Euler transformation to both sides of (2.19),
we obtain
\begin{multline}
(x\partial +b)x^{1-c}(1-x)^{c-a-b}F(1-a, 1-b, 2-c;x)\\
=(b+1-c)x^{1-c}(1-x)^{c-a-b-1}F(1-a, -b, 2-c;x),
\end{multline} 
and thus, 
\begin{multline}
H(\ell)x^{1-c}(1-x)^{c-a-b}F(1-a, 1-b, 2-c;x)\\
=(b+1-c,\ell)x^{1-c}(1-x)^{c-a-b-\ell}F(1-a, 1-b-\ell, 2-c;x).
\end{multline}
Substituting $b=1$ into (2.21), we obtain (2.18).
\section{A proof of the main theorem}
In this section, we prove Theorem 1.1 using Lemmas 2.2 and 2.3.

Operating with $H_1(\ell)$ on $F(a, 1, c; x)$,
from Lemma 2.3, we obtain
\begin{flalign}
&F(a, 1+\ell , c; x)=\frac{H_1(\ell)}{(1,\ell)}F(a, 1, c;x) \notag \\
&=\frac{1}{(1,\ell)}(q(x)\partial +r(x))\left( y_2(x)\int_0^{x}\frac{-(1-c)}{t(1-t)y_2(t)} dt\right) \notag \\ 
&=\frac{1}{(1,\ell)}\left[\left(q(x)\partial +r(x)\right)y_2(x)\times \int_0^{x}\frac{-(1-c)}{t(1-t)y_2(t)} dt
-\frac{(1-c)q(x)}{x(1-x)}\right]\notag \\
&=\frac{1}{(1,\ell)}\left[H_1(\ell)y_2(x)\times \int_0^{x}\frac{-(1-c)}{t(1-t)y_2(t)} dt
-\frac{(1-c)q(x)}{x(1-x)}\right]\notag \\
&=\frac{(2-c,\ell )}{(1,\ell)}y_2(x)(1-x)^{-\ell}F(1-a, -\ell ,2-c;x) \int_0^{x}\frac{-(1-c)}{t(1-t)y_2(t)} dt
-\frac{(1-c)q(x)}{(1,\ell)x(1-x)}
\end{flalign}
for $\Re c>1$.
In particular, for any root $\la$ of the polynomial $F(1-a, -\ell , 2-c; x)$, we have
\begin{gather}
F(a, 1+\ell , c; \la)=-\frac{(1-c)q(\la)}{(1,\ell)\la(1-\la)},
\end{gather}
and thus, from (2.12), this implies
\begin{gather}
F(a, 1+\ell , c; \la)=-\frac{(1-c)q_0(\la)}{(1,\ell)(1-\la)^{\ell}}.
\end{gather}
In addition, substituting $\la$ for $x$ in (2.13), we obtain
\begin{gather}
F(c-a, c-1-\ell , c; \la)=-\frac{(1-c)}{(1,\ell)}(1-\la)^{a+1-c}q_0(\la).
\end{gather}
Here, we have assumed that $\Re c>1$.
However, because both sides of (3.3) are analytic functions of $c$,
(3.3) is valid even if this condition is not satisfied,
by virtue of analytic continuation. 
The same holds for (3.4).
This completes the proof of Theorem 1.1.

Finally, using a symmetry of $q_0(x)$, we confirm the statement given in Remark 1.2.
We begin by assuming that $a\notin \Z$.
Then, although we have expressed $q_0(x)$ as in (2.13),
it can also be expressed as 
\begin{multline}
q_0(x)=x^{\ell -1}\times \\
\left[
-\frac{(2-c,\ell )(1,\ell )}{(a-1)(2-a,\ell)}\left(\frac{1}{x}\right)^{\ell}
\left(1-\frac{1}{x}\right)^{c-a-1}F\left(c-a, 1-a, 2-a+\ell; \frac{1}{x}\right)\right. \\
\left. +\frac{(-1)^{\ell}(1-a,\ell)}{a-1}F\left(2-c, 1, 2-a; \frac{1}{x}\right)
F\left(c-1-\ell , -\ell , a-\ell ;\frac{1}{x}\right)\right]
\end{multline}
(cf. Theorems 3.7 and 3.8 of [Eb]).
Note that the right-hand side of (3.5) is meaningful by the above assumption.
Therefore, we obtain (1.7). 

\textbf{Acknowledgement.}
The author would like to thank Professor Hiroyuki Ochiai and Professor Masaaki Yoshida
for many valuable comments and Professor Raimundas Vidunas for 
providing the author with a series of letters written by R.W.Gosper to D. Stanton.
The author also thanks the editor and referees for their constructive suggestions and comments.

\medskip
\begin{flushleft}
Akihito Ebisu\\
Department of Mathematics\\
Kyushu University\\
Nishi-ku, Fukuoka 819-0395\\
Japan\\
a-ebisu@math.kyushu-u.ac.jp
\end{flushleft}

\end{document}